\documentclass[11pt]{article}
\usepackage{amsmath,amssymb,theorem,enumerate,color}
\usepackage{graphicx}
\usepackage{float,color,fancybox,shapepar,setspace,hyperref}
\usepackage[affil-it]{authblk}
\usepackage{caption}
\usepackage{subfigure}
\usepackage{psfrag}
\usepackage{epstopdf}
\newtheorem{thm}{Theorem}[section]
\newtheorem{conj}[thm]{Conjecture}
\newtheorem{cor}[thm]{Corollary}
\newtheorem{lem}[thm]{Lemma}
\theorembodyfont{\rmfamily}
\def\pf{\medskip\noindent {\bf Proof.}~~}

\def\dfn#1{{\sl #1}}

\newcommand{\less}{\setminus}
\def\qed{ \hfill $\blacksquare$}
\def\es{\emptyset}
\newcounter{counter}

\begin{document}

\title{Gallai-Ramsey number  of even cycles with chords}

\author{Fangfang  Zhang$^{1,}$\thanks{This work was done while the first author visited the University of Central Florida as a visiting student. The hospitality of the
hosting institution is greatly acknowledged.    The visit was funded by the Chinese Scholarship Council.  E-mail address: fangfangzh@smail.nju.edu.cn.}  \,,  Zi-Xia Song$^{2,}$\thanks{Partially supported by the National   Science  Foundation of China under Grant No. DMS-1854903.   E-mail address:  Zixia.Song@ucf.edu. }\,\,,    Yaojun Chen$^{1,}$\thanks{Yaojun Chen and Fangfang Zhang are partially supported by the National Natural Science Foundation
of China  under grant numbers 11671198 and 11871270. E-mail address:  yaojunc@nju.edu.cn.}}
 \affil{ { \small {$^1$Department of Mathematics, Nanjing University, Nanjing 210093, P.R. China}}\\
  { \small {$^2$Department  of Mathematics, University of Central Florida, Orlando, FL 32816, USA}}
 }
 \date{ }
 \maketitle

\begin{abstract}
       For a graph $H$ and an integer $k\ge1$, the $k$-color Ramsey number $R_k(H)$ is   the least integer $N$ such that every $k$-coloring of the edges of the complete graph $K_N$  contains a monochromatic copy of $H$.   Let $C_m$ denote the cycle  on  $m\ge4$ vertices        and    let $\Theta_m$ denote  the family of   graphs obtained from   $C_m$ by adding  an additional edge joining two non-consecutive vertices.    Unlike Ramsey number of odd cycles,  little is known about the general behavior of $R_k(C_{2n})$  except  that $R_k(C_{2n})\ge (n-1)k+n+k-1$ for all $k\ge2$ and $n\ge2$.    In this paper,  we  study Ramsey number of even cycles with chords under Gallai colorings, where a  Gallai coloring  is a coloring of the edges of a complete graph without rainbow triangles.     For an integer $k\geq 1$,   the  Gallai-Ramsey number $GR_k(H)$ of a graph $H$ is the least positive  integer $N$ such that every Gallai $k$-coloring of the complete graph $K_N$   contains a monochromatic copy of $H$.  We prove that  $GR_k(\Theta_{2n})=(n-1)k+n+1$ for all $k\geq 2$ and $n\geq 3$. This implies that 
 $GR_k(C_{2n})=(n-1)k+n+1$  all $k\geq 2$ and $n\geq 3$.   Our result yields     a unified proof for the Gallai-Ramsey number of  all even cycles on at least four vertices.  
 \medskip
 

\end{abstract}
\baselineskip 16pt
\section{Introduction}
 In this paper we consider graphs that are finite, simple and undirected.      We use $P_m$,  $C_m$ and $K_m$ to denote the path,    cycle and  complete graph  on $m$ vertices, respectively.  For $m\ge4$, let $\Theta_m$ denote  the family of   graphs obtained from   $C_m$ by adding  an additional edge joining two non-consecutive vertices.   
For any positive integer $k$, we write  $[k]$ for the set $\{1, \ldots, k\}$. 
 Given an integer $k \ge 1$ and graphs $H_1,  \ldots, H_k$, the classical \dfn{Ramsey number} $R(H_1,   \ldots, H_k)$   is  the least    integer $N$ such that every $k$-coloring of  the edges of      $K_N$  contains  a monochromatic copy of  $H_i$ in color $i$ for some $i \in [k]$.  When $H = H_1 = \dots = H_k$, we simply write   $R_k(H)$ to denote the \dfn{$k$-color Ramsey number} of $H$.    In the seminal paper of Ramsey \cite{Ramsey}, it is shown that Ramsey numbers are finite. This was   rediscovered independently by Erd\H{o}s and Szekeres  \cite{ES}.  Since the 1970s, Ramsey theory has grown into one of the most  active areas of research   in combinatorics, overlapping variously with graph theory, number theory, geometry and logic. \medskip
 
Let $(G,\tau)$ denote a $k$-colored  complete graph, where $G$ is a complete graph and $\tau: E(G)\rightarrow [k]$. We say $(G,\tau)$ is \dfn{$\mathcal{F}$-free} if $G$ does not contain a monochromatic copy of a   graph in a given family $\mathcal{F}$ under  the $k$-coloring $\tau$; when $\mathcal{F}=\{F\}$, we simply say $(G,\tau)$ is \dfn{$F$-free}. By abusing notation, we say  $(G,\tau)$ contains a \dfn{monochromatic copy of   $\mathcal{F}$} if $G$  contains a monochromatic copy of a   graph in  $\mathcal{F}$ under $\tau$.
  One of the earliest and well-known problems is that of  determining the Ramsey number    $R_k(C_m)$. 
When $k=2$, the behavior of the Ramsey number $R(C_\ell,C_m)$ has been studied and fully determined by several authors, including Bondy and Erd\H{o}s \cite{BE}, Faudree and Schelp \cite{FS}  and Rosta \cite{Rosta}.   However, in the case where more than two colors are involved such results are still rather rare.   
 For even cycles, not much is known about the behavior of $R_k(C_{2n})$ in general.  Figaj and \L uczak \cite{FL} proved that for $\alpha_1,\alpha_2,\alpha_3 >0$,
  $$R(C_{2\lfloor\alpha_1n\rfloor},C_{2\lfloor\alpha_3n\rfloor},C_{2\lfloor\alpha_3n\rfloor})=(\alpha_1+\alpha_2+\alpha_3 +\max\{\alpha_1,\alpha_2,\alpha_3 \})n+o(n),$$
  as $n\rightarrow \infty$. 
  Following the ideas of   Gy\'arf\'as, Ruszink\'o, S\'ark\"ozy, and Szemer\'edi~\cite{GRSS} in determining the value of $R(P_m, P_m, P_m)$,  Benevides and  Skokan \cite{BS}   proved  that 
$R(C_{2n},C_{2n},C_{2n})=4n $ for sufficiently large $n$.
For general $k$, Dzido, Nowik and Szuca~\cite{DNS} showed that 

 \begin{thm}[\cite{DNS}]\label{EvenCycle} For all  $k\ge1$ and $n\ge2$,

\[
R_k(C_{2n }) \ge \begin{cases}
			(n-1)k+n+k  & \text{if } k \text{ is odd} \\
			(n-1)k+n+k-1 & \text{if } k \text{ is even.}
			\end{cases}
\]

 \end{thm}
For further results in this direction, we refer the reader to Graham, Rothchild and Spencer~\cite{GRS} and Radziszowski    \cite{survey} for  a dynamic survey.\medskip

In this paper we investigate  Ramsey numbers of even cycles and even cycles with chords under Gallai colorings,   
 where a \dfn{Gallai coloring} is a coloring of the edges of a complete graph without rainbow triangles (that is, a triangle with all its edges colored differently). Gallai colorings naturally arise in several areas including: information theory~\cite{KG}; the study of partially ordered sets, as in Gallai's original paper~\cite{Gallai} (his result   was restated in \cite{Gy} in the terminology of graphs); and the study of perfect graphs~\cite{CEL}. There are now a variety of papers  which consider Ramsey-type problems in Gallai colorings (see, e.g., \cite{C9C11, C13C15, C7, chen,  c5c6,GS, exponential, Hall,   C10C12, K4,   C6C8,W4, oddchorded}).     More information on this topic  can be found in~\cite{FGP, FMO}.  \medskip
 
A \dfn{Gallai $k$-coloring} is a Gallai coloring that uses at most $k$ colors. 
 Given an integer $k \ge 1$ and graphs $H_1,  \ldots, H_k$,  \dfn{Gallai-Ramsey number}  $GR(H_1,  \ldots, H_k)$ is defined to be  the least integer $N$ such that every Gallai $k$-coloring of $K_N$   contains a monochromatic copy of $H_i$ in color $i$ for some $i \in [k]$. When $H = H_1 = \dots = H_k$, we simply write $GR_k(H)$.    Clearly, $GR_k(H) \leq R_k(H)$ for all $k\ge1$ and $GR(H_1, H_2) = R(H_1, H_2)$.     
Theorem~\ref{general} below is a result of Gy\'{a}rf\'{a}s,   S\'{a}rk\"{o}zy,  Seb\H{o} and   Selkow~\cite{exponential} which characterizes    the general behavior of $GR_k(H)$.

\begin{thm} [\cite{exponential}]\label{general}
Let $H$ be a fixed graph  with no isolated vertices 
 and let $k\ge1$ be an integer. Then
$GR_k (H) $ is exponential in $k$ if  $H$ is not bipartite,    linear in $k$ if $H$ is bipartite but  not a star, and constant (does not depend on $k$) when $H$ is a star.				
\end{thm}

It turns out that for some graphs $H$ (e.g., when $H=K_3$),  $GR_k(H)$ behaves nicely, while the order of magnitude  of $R_k(H)$ seems hopelessly difficult to determine.   
We will utilize the following important structural result of Gallai~\cite{Gallai}.

\begin{thm}[\cite{Gallai}]\label{Gallai}
	Let  $(G, \tau) $  be a Gallai $k$-colored complete graph    with   $|V(G)|\ge  2 $.  Then  $V(G)$ can be partitioned into nonempty sets  $V_1,  \dots, V_p$ with $p\ge2$ so that    at most two colors are used on the edges in $E(G)\less (E(V_1)\cup \cdots\cup  E(V_p))$ and only one color is used on the edges between any fixed pair $(V_i, V_j)$ under $\tau$, where $E(V_i)$ denotes the set of edges with both ends in  $ V_i $ for all $i\in [p]$. 
\end{thm}

 The partition $\{V_1,   \dots, V_p\}$ given in Theorem~\ref{Gallai} is  a \dfn{Gallai partition} of    $(G, \tau)$.  Let $(\mathcal{R},\tau)$ be obtained from  $(G, \tau)$ by first contracting each $V_i$ into a single vertex $v_i$ and then coloring $v_iv_j$ by the unique color on the edges between $V_i$ and $V_j$ in $(G, \tau)$.   We say  $(\mathcal{R},\tau)$ is   the  \dfn{reduced graph} of $(G,\tau)$ corresponding to the  Gallai partition $\{V_1,   \dots, V_p\}$. Note that   $\mathcal{R}=K_p$.  
By Theorem~\ref{Gallai},  all edges in $\mathcal{R}$ are colored by at most two colors under $\tau$.  One can see that any monochromatic copy of $H$ in $(\mathcal{R},\tau)$  will result in a monochromatic copy of $H$ in $(G, \tau)$. It is not  surprising  that    the  $2$-color  Ramsey number  $R_2(H)$ plays an important role in determining the value of  $GR_k( H)$ when $H$ is a complete graph.    Fox,  Grinshpun and  Pach~\cite{FGP} posed the following  conjecture.  

\begin{conj}[\cite{FGP}]\label{Fox} For all  $k\ge1$ and $t\ge3$,
\[
GR_k( K_t) = \begin{cases}
			(R_2(K_t)-1)^{k/2} + 1 & \text{if } k \text{ is even} \\
			(t-1)  (R_2(K_t)-1)^{(k-1)/2} + 1 & \text{if } k \text{ is odd.}
			\end{cases}
\]
\end{conj}

 The first case of   Conjecture~\ref{Fox} follows directly from a result of Chung and Graham \cite{chgr} in 1983. The proof in \cite{chgr}  does not rely on Theorem~\ref{Gallai}.  
 A simpler proof of the case $t=3$ using Theorem~\ref{Gallai}  can be found in~\cite{exponential}.    The next open case, when $t=4$,  was recently settled in~\cite{K4}.    \medskip

 For the remainder of the paper, 
 we use $(G, \tau)$ to denote a Gallai $k$-colored complete graph, where $G$ is a complete graph and   $\tau: E(G)\rightarrow [k]$ is a Gallai $k$-coloring.   For  each $(G, \tau)$, let  $G^\tau_i$   denote the spanning subgraph of $G$ with   $E(G^\tau_i): = \{e\in E(G)\mid  \tau(e)=i\}$ for all $i\in[k]$. We simply write $G^\tau_r$ if the color $i$ is red; $G^\tau_b$ if the color $i$ is blue.   
For  every   $n\in \mathbb{N}$,    let  $q_\tau(G, n)$ denote   the number of colors $i\in[k]$ such that  $G^\tau_i$ has  a component of order at
least $n$. Then  $q_\tau(G, n)\le k$.   
We begin with  Lemma~\ref{k'} ( which follows directly from the proof of Lemma 9 given in \cite{Hall}), which we restate and prove here for completeness.  

\begin{lem}[\cite{Hall}]\label{k'}
Let  $(G, \tau) $  be a Gallai $k$-colored complete graph    with   $|V(G)|\ge n\ge2 $. Then the following hold.  \vspace{-0.15cm}
\begin{enumerate}[\rm(i)] 
\item\label{qtau1}   $q_\tau(G, n)\ge1$.  
\vspace{-0.15cm}
 \item\label{qtau2}  Let $\{V_1, \ldots, V_p\}$ be  a Gallai partition of $(G, \tau)$   with  $p \ge 2$ as small as possible. Then for any color $i$ on the edges  in $(\mathcal{R}, \tau)$, $\mathcal{R}^\tau_i$   is connected.   
 \end{enumerate}
\end{lem}
\vspace{-0.2cm}
\pf    Let     $\{V_1, \ldots, V_p\}$ be  any Gallai partition of $(G,\tau)$ and let $(\mathcal{R},\tau)$ be the corresponding reduced graph   of $(G,\tau)$.  By Theorem~\ref{Gallai}, we may further assume that every edge  of $\mathcal{R}$ is  colored red or blue. Note that for any $i\in\{r,b\}$, $G_i^{\tau}$ is connected if  $\mathcal{R}^{\tau}_i$ is connected.   Hence the statement is trivially true if $p=2$ or  $(\mathcal{R},\tau)$ is monochromatic. We may then assume that $p\ge3$ and     $(\mathcal{R},\tau)$  contains  both red and  blue edges.      If  both   $\mathcal{R}^{\tau}_r$ and $\mathcal{R}^{\tau}_b$  are connected, then we are done, so we may assume that  $\mathcal{R}^{\tau}_r$ is disconnected. Let $D$ be a component of  $\mathcal{R}^{\tau}_r$.   Then all edges between $V(D)$ and $V(\mathcal{R})\less V(D)$ in $(\mathcal{R}, \tau)$ are blue, and so  $\mathcal{R}^{\tau}_b$ must be connected. Hence $q_\tau(G, n)\ge1$.  Next assume that $p\ge3$ is chosen as small as possible.  We may assume that $V(D)=\{v_1, \ldots, v_{d}\}$, where $d<p$ and $v_1, \ldots, v_d$ are obtained from $V_1, \ldots, V_d$ by contracting each $V_i$ into $v_i$.  Let  $V_1':=\cup_{j=1}^d V_j$ and $V_2':=V(G)\less V_1'$. Then $\{V_1', V_2'\}$ yields a Gallai partition of $(G,\tau)$ with two parts, contrary to the minimality of $p$. Consequently, if  $p $ is chosen as small as possible, then  $\mathcal{R}^\tau_i$   is connected for any color $i$ on the edges  in $(\mathcal{R}, \tau)$. \qed\medskip

In this paper, we   focus on determining the exact values  of  $GR_k(\Theta_{2n})$ for all $k\ge2$ and $n\ge3$, and $GR_k(C_{2n})$ for all $k\ge2$ and $n\ge2$.  Note that $GR_k(\Theta_{2n})\ge GR_k(C_{2n})$ for all $k\ge2$ and $n\ge2$.  Using the construction of    Erd\H{o}s, Faudree, Rousseau and Schelp    (see Section 2  in \cite{EFRS}) for $R_k(C_{2n})$, we have   \[GR_k(C_{2n })\ge\begin{cases}
			 (n-1)  k+n+1    & \text{if }   k\ge2 \text{ and } n\ge3  \\
			   (n-1) k+n+2    & \text{if }  k\ge2 \text{ and } n=2.
\end{cases} \]
 Indeed when $k=2$, simply take $(G,\tau)$ to be a 2-colored $K_5$ with each color being a  monochromatic $C_5$ when $n=2$, and a 2-colored $K_{3n-2}$ obtained from 1-colored  $K_{2n-1}$  by first adding $n-1$ vertices  and then coloring  all the new edges with  the second  color  when $n\ge3$;  for $k\ge 3$, let $(G,\tau)$ be obtained from  the construction for $k -1$ by first adding $n-1$ vertices and then coloring all the new edges with a new color.      General upper bounds for $GR_k(C_{2n})$ were first studied in \cite{c5c6}, later improved in \cite{Hall}.
 More recently,   the exact value of $GR_k(C_{2n }) $ for $2\le n\le6$    has been completely settled, see  \cite{FGJM, c5c6,  c8, C10C12, C6C8}. Building on the ideas in \cite{Hall}, we establishes  a sufficient condition for the existence of monochromatic copy of  $\Theta_{2n}$ in   Gallai $k$-colored complete graphs. We include $C_4$ in the statement of Theorem~\ref{even} below  in order to   provide     a unified proof for the Gallai-Ramsey number  of all even cycles on at least four vertices.

\begin{thm}\label{even}
 Let  $(G, \tau) $  be a Gallai $k$-colored complete graph    with   $|V(G)|\ge n\ge2 $. If  
\[|V(G)| \ge \begin{cases}
			 (n-1)\cdot q_\tau(G,  n)+n+1    & \text{if }   n\ge3  \\
			   (n-1) \cdot q_\tau(G,  n)+n+2    & \text{if }  n=2, 
			\end{cases} \]  
then $(G, \tau)$ has a  monochromatic   $C_4$ when $n=2$ and a  monochromatic  $\Theta_{2n}$ for all  $n\ge3$.
\end{thm}

 Theorem~\ref{even} turns out to be very powerful. Since   $q_\tau(G, n)\le k$, we see that Theorem~\ref{even}   yields the exact values of $GR_k(\Theta_{2n})$ for all $k\ge2$ and $n\ge3$;  exact values of  $GR_k(C_{2n})$    for all $k\ge2$ and $n\ge 2$.  Furthermore, Theorem~\ref{even} also provides  a unified  proof for the Gallai-Ramsey number  of $C_{2n}$ for all $k\ge2$ and $n\ge2$.

 \begin{thm}\label{evenchorded}
For   all     $k\geq 2$ and $n\geq 3$,   
 we have $GR_k(\Theta_{2n })=  (n-1)  k+n+1$.    
\end{thm}

 \begin{thm}\label{evenexact}
For   all     $k\geq 2$ and $n\geq 2$,    we have 
\[GR_k(C_{2n })=\begin{cases}
			 (n-1)  k+n+1    & \text{if }   n\ge3  \\
			   (n-1) k+n+2    & \text{if }  n=2.
\end{cases} \]
\end{thm}

 Using a similar argument in   \cite[Proposition 1.14]{C6C8},  Theorem~\ref{evenexact} further yields  the exact value of  $GR_k(P_m)$  for all $k\ge2$ and $m\ge6$ (see  \cite[Proposition 1.12]{C6C8} for the lower bound construction for $GR_k(P_m)$). 

 \begin{thm}\label{pathexact}
For   all     $k\geq 2$ and $m\geq 6$,   
\[GR_k(P_m)=\begin{cases}
			 (n-1)  k+n+1    & \text{if }   m=2n   \\
			   (n-1) k+n+2    & \text{if }  m=2n+1.
\end{cases} \]
\end{thm}

We prove Theorem~\ref{even} in Section~\ref{GR-even}. The proof of Theorem~\ref{even} uses recoloring method and several structural results  on    the existence of  a cycle or path in bipartite graphs. An overview of the proof is given in Section~\ref{overview}.        \medskip

We conclude this section by introducing more notation.   Given a graph $G$,    sets $S\subseteq V(G)$ and $F\subseteq E(G)$,  we use   $|G|$    to denote  the  number
of vertices    of $G$,      $G\less S$ the subgraph    obtained from $G$ by deleting all vertices in $S$,   $G\less F$ the subgraph    obtained from $G$ by deleting all edges in $F$, and $G[S] $    the  subgraph    obtained from $G$ by deleting all vertices in $V(G)\less S$.  We simply write $G\less v$ when $S=\{v\}$, and $G\less uv$ when $F=\{uv\}$.   
      For  two disjoint sets $A, B\subseteq V(G)$,    $A$ is \dfn{complete} to $B$ in $G$  if each vertex in $A$ is adjacent to all vertices in  $B$, and \dfn{anti-complete} to $B$ in $G$  if no  vertex in $A$ is adjacent to any vertex  in  $B$.   
 Let $(G,\tau)$ be a Gallai $k$-colored complete graph.   
For  two disjoint sets $A, B\subseteq V(G)$,   $A$ is \dfn{mc-complete} to $B$   if all the edges between $A$ and $B$  in    $(G,\tau)$ are colored the same color. We simply say      $A$ is     \dfn{$j$-complete} to $B$    if all the edges between $A$ and $B$  in   $(G,\tau)$ are colored by some color $j\in[k]$,   and   $A$ is \dfn{blue-complete}     to $B$    if all the edges between $A$ and $B$  in $(G,\tau)$ are colored  blue.  We say a vertex $x\in V(G)$ is \dfn{blue-adjacent} to a vertex $y\in V(G)$ if the edge $xy$ is colored blue in $(G,\tau)$, and $x$ is \dfn{blue-complete} to an edge $yz\in E(G)$ if $x$ is  blue-complete to $\{y, z\}$ in $(G,\tau)$.  
Similar definitions hold when blue is replaced by another color. For convenience, we use  $A \less B$ to denote  $A-B$; and  $A \less b$ to denote  $A \less \{b\}$ when $B=\{b\}$. We use the convention   ``$S:=$'' to mean that $S$ is defined to be the right-hand side of the relation.

  \section {An overview of the proof of Theorem~\ref{even}}\label{overview}

Let $(G,\tau)$, $n$ and $k$ be as given  in the statement of Theorem~\ref{even}. Suppose the statement is false.  Choose $n\ge 2$ as small as possible, and subject to the choice of $n$,  choose $(G,\tau)$ so that $q_\tau(G, n)$  is minimum.   Let $q:=q_\tau(G, n)\ge1$.  We may assume that for each color $i\in [q]$, $G^\tau_i$ has  a component of order at
least $n$.    Let $X_1,\ldots, X_q$ be disjoint subsets of $V(G)$ such that for each $i\in [q]$, $X_i$ (possible empty) is mc-complete in color $i$ to $V(G)\less \bigcup_{i=1}^qX_i$. Choose $X_1,\ldots, X_q$ so that $|G|-\sum_{i=1}^q  |X_i|\geq n$ and $\sum_{i=1}^q  |X_i|$ is as large as possible. Let $X:=\bigcup_{i=1}^qX_i$ and let $\{V_1, \ldots, V_p\}$ be  a Gallai partition of $(G \less X, \tau)$ with       $p \ge 2$ as small as possible and   $|V_1|  \le \cdots \le |V_p|$.   Then    $(G,\tau)$ is   $C_{2n}$-free and    the reduced graph $(\mathcal{R}, \tau)$ of 
$(G \less X, \tau)$ are colored by at most two colors in $[q]$, say red and blue. Let $R$ be set of all vertices   $v\in V(G)\less (X\cup V_p)$ such that $v$ is red-complete to $V_p$, and $B$ be set of all vertices   $v\in V(G)\less (X\cup V_p)$ such that $v$ is blue-complete to $V_p$. With some effort it can be shown
that $R\ne \es$,  $B\ne \es$, $q=2$,  $n\ge3$ and $|X_i|\le n-1$ for all $i\in[2]$. We may assume $X_1$ is red-complete to $V(G)\less X$ and $X_2$ is blue-complete to $V(G)\less X$. We then prove  two crucial claims  (5) and (7)  that $  |V_p|\le n-2$ and either $X_1=\es$ or $X_2=\es$. These allow us to reduce $\tau$ to be a Gallai $3$-coloring as follows:  
let $\sigma$ be  obtained from $\tau$ by recoloring all the edges of $ G[X_1] $ blue if $X=X_1$ and $ G[X_2] $ red if $X=X_2$; and    all the edges of  $ G[V_i] $ green for all $i\in[p]$. Then $(G, \sigma)$ is   Gallai $3$-colored   with no monochromatic copy of $C_{2n}$,    $q_\sigma(G, n)=q=2$ and   $\{V_1, \ldots, V_p\}$  a Gallai partition of  $(G\less X, \sigma)$. Let $x\in V_p$, $y\in R \cup X_1$ and $z\in B\cup X_2 $ be  such that $y\in X_1$ if $X_1 \ne\es$ and $z\in X_2$ if $X_2 \ne\es$. Let $H:=G\less \{x,y,z\}$.
 Then $q_\sigma(H, n-1)\leq q_\tau(H, n-1) =2$. By minimality of $n$,  $(H, \sigma)$ contains a red or blue copy of $ C:=C_{2n-2}$, say blue. We then define  three pairwise disjoint  sets $S, T, L$ according to the green edges between $V(C)$ and $V(G)\less V(C)$ in $(G, \sigma)$:  $S$ consists of all (special) vertices $v\in V(C)$  such that $v$ is green-adjacent to some vertex  in $V(G)\less V(C)$; $T$  consists of   all vertices  $u\in   (R \cup B)\less V(C)$ such that     $u$ is green-adjacent to   some vertex in $  S$; finally, let $L:=V(G)\less (V(C) \cup T)$, all leftover vertices.  In the next step we prove  that $|L|\le n$,  $|T|\ge1$ and $1\le |S|\le n-4$ by either finding a red $C_{2n}$ using edges between $V(C)$ and $T\cup L$ or extending $C$ to be a blue $C_{2n}$.  The  key claim (15)  states that:  there exists a subgraph $J$ of $C\less S$ with $|J|\ge 2n-2-2|S|\ge6$ such that    each component of $J$ is a path of odd length, and   for all $u\in   T\cup L$ except possiblely one vertex,    $u$ is red-adjacent to at least $|J|/2$ many    vertices in $J$.  Finally, in the last step,  we choose ``wisely" a subset $W$ of $T\cup L$  and  a subset $W^*$ of $(T\cup L)\less W$  with  $|W|=n-1-s$ and $|W^*|=s$.  With some efforts  it can be shown that  $(G[W^*\cup S],\sigma)$  has a red cycle $C^*$ on $2s$ vertices, and $(G[W \cup V(J)],\sigma)$ has a red path $P$ on $2n-(2s+3)$ vertices with both   ends in $W$ such that $P$ uses only edges   between $W$ and $V(J)$. By the ``wise" choice of $W$ and $W^*$, we obtain a red $C_{2n}$ from  $C^*$ and $P$ by joining them through three additional vertices.

\section{ Preliminaries}

  In this section we first list some known results and then prove several new corollaries that shall be applied in the proof of our main results.

\begin{thm}[\cite{Rosta}]\label{C2n}
For all   $n\ge3$, $   R(C_{2n}, C_{2n}) = 3n-1$.
\end{thm}
\begin{thm} [\cite{Bon}]\label{Bon}
Let $G$ be a graph  on $n\ge3$ vertices. If $\delta(G)\geq n/{2}$, then either $G$ has a cycle of length $\ell$ for every $\ell$ satisfying $3\le \ell\le n$ or $n$ is even and $G$ is isomorphic to $K_{ n/2, n/2 }$.
\end{thm}

 For a bipartite graph $G$ with bipartition $\{M, N\}$, let     $\delta(M):= \min\{d_G(x): x\in M\}$ and $\Delta(M):= \max\{d_G(x): x\in M\}$.  Lemma~\ref{Hal3}   follows from Lemma 12, Lemma 13, Lemma 14 in \cite{Hall}. 
 
\begin{lem}[\cite{Hall}]\label{Hal3}
Let $G$ be a bipartite graph  with bipartition $\{M, N\}$ such that $|M| \geq 2$, $|N| \geq 4$ and  $\delta(M) \ge (|N| + 1)/{2}$. Then $G$ has a cycle of length $ 2\ell $ for any $\ell$ satisfying $2 \leq \ell \leq \min\{|M|, \delta(M)-1\}$,  or $\Delta(M)=(|N|+1)/2$, $M=M_1\cup M_2$ and  $ N=N_1\cup  N_2\cup N_3$, where $M_1, M_2, N_1,  N_2, N_3$ are non-empty, pairwise disjoint  sets in $G$, $|N_3|=1$, $|N_1|=|N_2|$, and $M_i  $ is  complete to $N_i\cup N_3$ but anti-complete to $N_{3-i}$ for all $i\in[2]$.    \end{lem} 

We next prove  several structural results on      the existence of  a cycle or path in bipartite graphs,   which will be useful in the proof of  our main results. 

\begin{cor}\label{P2l-13}
Let $G$ be a bipartite graph  with bipartition $\{M, N\}$ such that $|M| \geq 2$, $|N| \geq3$  and  $\delta(M) \ge |N|/{2}$. Let $\ell = \min\{|M|, \delta(M)\}$. 
   Then $G$ has a path on $2\ell-1$  vertices with both ends in $M$ or $M=M_1\cup M_2$ and  $ N=N_1\cup  N_2$, where $M_1, M_2, N_1,  N_2 $ are non-empty, pairwise disjoint sets, $|N_1|=|N_2|$, and $M_i  $ is  complete to $N_i $ but anti-complete to $N_{3-i}$ for all $i\in[2]$.    Moreover, if    $|M|\ge 2\ell-3$ and  $|N|=2\ell$, then    $G$ has a path on $2\ell -3$ vertices with both ends in $N$. 

\end{cor}

\pf Let $H$ be obtained from $G$ by adding a new vertex $x$ adjacent to all vertices in $M$. Then $H$ is a bipartite graph with bipartition $\{M, N\cup \{x\}\}$. For any vertex  $v\in M$,   $d_H(v)\ge d_G(v)+1\ge (|N\cup \{x\}|+1)/2$.  By Lemma  \ref{Hal3}, $H$ contains a cycle $C_{2\ell}$ (and thus $G$ has a desired path on $2\ell-1$  vertices with both ends in $M$ and a desired path on $2\ell-3$  vertices with both ends in $N$) or  $M=M_1\cup M_2$ and $ N=N_1\cup  N_2$, where $M_1, M_2, N_1,  N_2 $ are non-empty, pairwise disjoint sets,   $|N_1|=|N_2|=\ell$, and $M_i  $ is  complete to $N_i $ but anti-complete to $N_{3-i}$ for all $i\in[2]$.    Moreover,  if      $|M|\ge 2\ell-3$ and  $|N|=2\ell$,  then $|M_1|+ |M_2|=|M|\ge 2\ell-3$, we see that $G$ contains a  path on $2\ell -3$ vertices with both ends in $N$, as desired.\qed

\begin{cor}\label{path}
Let $G$ be a bipartite graph  with bipartition $\{M, N\}$ such that $|N|\ge |M| \geq 3$ and  $\delta(M) \ge |N|-1$. If $M$ has at least two vertices each complete to $N$,   then $G$ contains a cycle of length  $2|M|$.
\end{cor}
\pf Let $G$, $M$ and $N$ be as in the statement. Let $M'$    be the set of all vertices $v\in M $  such that $v$  is complete to $N$. Then $|M'|\ge 2$.  Since $\delta(M)\ge |N|-1$, the statement is trivially true if $|M\less M'|\le  1$. So we may assume that   $|M\less M'|\ge 2$. By Corollary \ref{P2l-13} applied to $G\less M'$ and the fact $|N|-1>|N|/2$, $G\less M'$ has an $(x,y)$-path $P$ on $2|M\less M'|-1$ vertices with $x,y\in M\less M'$. Then $|N\less V(P)|\ge |M'|+1\ge 3$. Let $u,w \in N\less V(P)$ be distinct such that $xu, yw\in E(G)$ and let $Q$ be a $(u, w)$-path using edges between $M'$ and  $N\less V(P)$.  Then $G$ has a cycle on $2|M|$ vertices with edge set $E(P)\cup E(Q)\cup \{xu, yw\}$.\qed

\begin{cor}\label{C2l}
Let $G$ be a bipartite graph  with bipartition $\{M, N\}$ such that $|M|\geq 2 $ and $|N| \geq 4$.    If $\delta(M) > (|N| + 1)/{2}$  or  $\Delta(M) >  \delta(M) = (|N| + 1)/{2}$,  then $G$ has a cycle of length $ {2\ell}$  for any $\ell$ satisfying $2 \leq \ell \leq  \min\{|M|, \delta(M)\}$.  
\end{cor} 

\pf Let $\delta:=\delta(M)$ and let $x$ be a vertex in $M$ with $d (x)=\Delta(M)$. Assume first that $\delta>|M|\ge2$. By our assumption and Lemma \ref{Hal3},  $G $ has a  cycle of length $ 2\ell $  for any $\ell$ satisfying $2 \leq \ell \leq  |M|$. Assume next that $|M|\ge \delta$. 
By assumption, $|M|\ge\delta\ge3$ and  $|M\less x|\ge \delta-1\ge2$. By Lemma \ref{Hal3},  $G $ has a  cycle of length $ 2\ell $  for any $\ell$ satisfying $2 \leq \ell \leq  \delta -1$, and     
$G\less x$ has a cycle $C:=C_{2\delta-2}$, say  with vertices  $a_1, b_1, a_2, b_2, \ldots, a_{\delta-1}, b_{\delta-1}$ in order, where $a_1, \ldots, a_{\delta-1}\in M$ and $b_1, \ldots, b_{\delta-1}\in  N$.     We next show that $G$ has a cycle of length    $2\delta$. Let $\overline{N}:=N\less V(C)$.  By the choice of $x$ and $C$,  we see that $x\in M\less V(C)$, $N(x)\cap \overline{N}\ne \es$ and  
$|\overline{N}|  = |N|-(\delta-1)   <d(x)$. 
 Then   $N(x)\cap V(C)\neq \emptyset$.  Let \[M^* :=\{a_i \in\{a_1, \ldots, a_{\delta-1}\} \mid  b_i \in N(x),    \text{ where } i\in[\delta-1]\}.\] 
    Assume first that there exists a vertex, say $b\in N(x)\cap \overline{N}$, such that $b$ is adjacent to some $a_i\in M^*$. Then  $G$ contains  a cycle of length $2\delta$   with edge set $E(C\less a_ib_i)\cup \{a_ib, bx, xb_i\}$, as desired. So we may assume that $M^*$ is anti-complete to $N(x)\cap \overline{N}$.  Since 
   \[|M^*|+ |N(x)\cap \overline{N}|=d(x)>|\overline{N}|=|N(x)\cap \overline{N}|+ |\overline{N}\less N(x)|, \]
 we have  $|M^*|>  |\overline{N}\less N(x)|$.  Note that for each $a_i\in M^*$, 
   $|N(a_i)\cap \overline{N}|=d(a_i)-(\delta-1)\ge 1$. It follows that $  \overline{N}\less N(x)\ne \es$,  $|M^*|\ge 2$, and there must exist a vertex, say $b\in \overline{N}\less N(x)$, such that $b$ is complete to  $\{a_i,a_j\}$, where  $a_i,a_j\in M^*$ with $i<j$.  Then $G$ contains a cycle of length $2\delta$   with   edge set $E(C\less \{a_ib_i, a_jb_j\})\cup \{a_ib, ba_j, b_ix,  xb_j\}$, as desired. \medskip
   
    This completes the proof of Corollary \ref{C2l}
\qed

\section{Proof of Theorem~\ref{even}}\label{GR-even}

Let $(G,\tau)$, $n$ and $k$ be   as given in the statement. Suppose $(G,\tau)$ is  $C_4$-free when $n=2$ and  $\Theta_{2n}$-free when $n\ge3$.   Choose $n\ge 2$ as small as possible, and subject to the choice of $n$,  choose $(G,\tau)$ so that $q_\tau(G, n)$  is minimum.  Let $q:=q_\tau(G, n)$.  By Lemma~\ref{k'}(i), $q\ge1$. 
We may assume that for each color $i\in [q]$, $G^\tau_i$ has  a component of order at
least $n$.   Let $X_1,\ldots, X_q$ be disjoint subsets of $V(G)$ such that for each $i\in [q]$, $X_i$ (possiblely empty) is mc-complete in color $i$ to $V(G)\less \bigcup_{i=1}^qX_i$. Choose $X_1,\ldots, X_q$ so that $|G|-\sum_{i=1}^q  |X_i|\geq n$ and $\sum_{i=1}^q  |X_i|$ is as large as possible. Let $X:=\bigcup_{i=1}^qX_i$. Then $|G\less X|\geq n\ge2$. Since $(G, \tau)$ has no rainbow triangle, we see that for $i, j\in [q]$ with $i\ne j$,  each edge between $X_i$ and $X_j$ is colored $i$ or $j$. We next prove a series of claims. \\\
\setcounter{counter}{0}

\noindent {\refstepcounter{counter} \label{noC2n}  (\arabic{counter}) }
  $(G,\tau)$ is $C_{2n}$-free for all $n\ge2$.

\pf  Suppose $(G,\tau)$ contains a monochromatic copy of  $C:=C_{2n}$,  
 say    with vertices $x_1,x_2,\ldots,x_{2n}$ in order. Then $n\ge3$. We may assume that all edges of $C$ are colored blue.  Then  no chord of $C$ is colored blue   because  $(G,\tau)$ is  $\Theta_{2n}$-free. We may further assume that $x_1x_3$ is colored red. Then $x_1x_j$ is colored red for all $j\in \{3, 4, \ldots, 2n-1\}$ because $(G, \tau)$ has a no rainbow triangle. It follows that all chords of  $C$  are colored red.  Let $H$ be the graph with $V(H)=V(C)$ and $E(H)$ consisting of all chords of $C$. Then $H$ is the complement of $C_{2n}$. It can be easily checked that   $H$  contains a chorded $C_{2n}$ because $n\ge3$.   Thus $(G,\tau)$ contains a red copy of   $\Theta_{2n}$, a contradiction.\qed\\

\noindent {\refstepcounter{counter} \label{Xi}  (\arabic{counter}) }  
For all $i \in [q]$, $|X_i|  \le n-1$.
 
 \pf Suppose $|X_i| \geq n$ for some color $i \in [q]$, say blue. By the choice of $X_1,\ldots,X_q$, $X_i$ is blue-complete to $V (G) \less X$. It follows that $(G,\tau)$ has a blue $C_{2n}$ using edges between $X_i$ and $V (G) \less X$, contrary to  (\ref{noC2n}). 
 \qed\\

Let $\{V_1, \ldots, V_p\}$ be  a Gallai partition of $(G \less X, \tau)$ with       $p \ge 2$ as small as possible.  We may assume that  $|V_1|  \le \cdots \le |V_p|$.  By Theorem~\ref{Gallai} and Lemma~\ref{k'}(ii),  all   edges of the reduced graph of 
$(G \less X, \tau)$ are colored by at most two colors in $[q]$, say red and blue.   Then  for all $i\in[p-1]$,   $V_i$ is  either  red- or  blue-complete to $V_p$ in  $(G\less X,\tau)$.  Let 
\[
\begin{split}
 \mathcal{V}_r &:=    \{V_i  \in \{V_1, \ldots, V_{p-1}\} \mid V_i \text{ is   red-complete to } V_p \text{ in } (G\less X, \tau)\} \text{ and}\\
 \mathcal{V}_b  &:=    \{V_i  \in \{V_1, \ldots, V_{p-1}\} \mid V_i \text{ is   blue-complete to } V_p \text{ in } (G\less X, \tau)\}.
\end{split}
\]

Let $R=\bigcup_{V_j\in  \mathcal{V}_r }V_j$ and $B=\bigcup_{V_j\in  \mathcal{V}_b}V_j$. Then $R \cup B =V(G)\less (X\cup V_p)$,  and $R $ and $B $ are disjoint.   We may further assume that $X_1$ is red-complete to $V(G) \less X$ and $X_2$ is blue-complete to $V(G) \less X$. \\

\noindent {\refstepcounter{counter}\label{B}  (\arabic{counter}) }  
$R \ne\es$ and $B \ne\es$.   

\pf Suppose $R =\es$ or $B =\es$, say the latter.   Since $p\ge2$, we see that  $|R |\ge1$ and $R $ is red-complete to $V_p$.  Then $|V_p|\le n-1$, else, let  $X'_1:=X_1\cup R $ and $X'_i:=X_i$ for all $i\in\{2, \ldots, q\}$.  But then  
$ |X'_1|+\cdots+|X_q'| =|X \cup R  | >|X|$,
 contrary to  the  choice of $X_1, \ldots, X_q$.  Similarly, $|R |\le n-1$.  By (\ref{Xi}), $|X_i|\le n-1$ for all $i\in[q]$. If $|V_p\cup R \cup X_1|\leq 2n-1 $, then  \[|G|=|V_p\cup R\cup X_1|+|X_2|+\cdots +|X_q|\leq (2n-1)+(n-1)(q-1)=(n-1)q+ n,\]
  contrary to the assumption that $|G|\ge (n-1)q+n+1$. Thus    $|V_p\cup R \cup X_1|\geq 2n $.  Let $H$ be the subgraph of $G$ with $V(H)=V_p\cup R \cup X_1$ and $E(H)$ consisting of all red edges in $(G[V_p\cup R \cup X_1],  \tau)$.  Then $\delta(H)\ge |H|/2$. By Theorem~\ref{Bon},  $H$  has  a red $C_{2n}$,  which yields a red $C_{2n}$ in   $(G, \tau)$, contrary to  (\ref{noC2n}).  Thus $R \ne\es$ and $B \ne\es$.   \qed\\

  By (\ref{B}) and   Lemma~\ref{k'}(ii),    both $G_r^\tau\less X$ and $G_b^\tau\less X$ are connected. Thus $q\ge2$.  By minimality of $p$, $R $ is neither red- nor blue-complete to $B $ in $(G\less X, \tau)$. Thus $p\ge4$  and so $|G\less X|\ge p\ge4$. \\

\noindent {\refstepcounter{counter}\label{Vpn-1}  (\arabic{counter}) }  
$|V_p|\le n-1 $ and so $X_i \neq \emptyset$ for every color $i\in [q]$ that is neither red nor blue.

\pf Suppose $|V_p|\ge n$.  Then every vertex in $R  \cup B  $ is either  red- or  blue-complete to $V_p$.   Let $X'_1:=X_1\cup R $, $X'_2:=X_2\cup B  $, and $X'_i:=X_i$ for all $i\in\{3, \ldots, q\}$.  But then  
$ |X'_1|+\cdots+|X_q'| =|X \cup R   \cup B  | >|X|$,
 contrary to  the  choice of $X_1, \ldots, X_q$.   This proves that $|V_p|\le n-1$. Next, suppose there exists a color $j\in[q]$ such that $j$ is neither red nor blue but   $X_j= \emptyset $. Then   no edges between  pairs of $X_1, \ldots, X_q$ are colored by color $j$ in $(G,\tau)$. But then  $G^\tau_j$ has no component of order at least $n$,  because   $|V_\ell |\le  |V_p|\le n-1$   for all $\ell\in[p-1]$, and $|X_i|\le n-1$ for all $i\in [q]$ by (\ref{Xi}), a contradiction.  
  \qed\\

\noindent {\refstepcounter{counter}\label{q=2}  (\arabic{counter}) }  
$n\ge3$ and $q=2$.   

\pf Suppose first $n=2$.  By  (\ref{Vpn-1}) and   (\ref{Xi}), $|V_p|=1$ and $|X_i|\le1$ for all $i\in [q]$.   Since $(G, \tau)$ is $C_4$-free, we see that $q\ge2$. Then 
$|G\less \bigcup_{i=3}^q X_i|\ge q+4- (q-2)=6$. Thus   $(G,\tau)$ contains a red or blue $C_4$ because $R_2(C_4)=6$ \cite{CS}, contrary to  (\ref{noC2n}).   Suppose  next $n\ge 3$ and $q\ne 2$. Then $q\ge 3$  and  the color $q$ is neither red nor blue.      By (\ref{Vpn-1}),  $X_q\neq \emptyset$ and $|V_p|\le n-1$. Thus   $q_\tau(G\less X_q, n)=q_\tau(G,n)-1=q-1$.   By (\ref{Xi}),  $|X_q|\leq n-1$, and so \[ |G\less X_q|\geq (n-1)q+n+1-(n-1) =(n-1)\cdot q_\tau(G\less X_q,n)+n+1.\] By  minimality  of $q$ and the fact that $n\ge3$, $(G\less X_q,\tau)$ has a monochromatic copy of $C_{2n}$, contrary to  (\ref{noC2n}). This proves that  $q=2$. 
\qed\\

  By (\ref{q=2}), $n\ge3$ and $q=2$ and so  $|G|\geq 2(n-1)+n+1=3n-1$.  By Theorem~\ref{C2n},   $k\ge3$ and  edges of $(G,\tau)$ must be colored by at least three colors.    We may further assume that   the third color on $E(G)$ is green under  $\tau$.     
\\

\noindent {\refstepcounter{counter}\label{Vpn-2}  (\arabic{counter}) }  
$ |V_p|\leq n-2$.  

\pf    Suppose $|V_p|\geq n-1$.  By  (\ref{Vpn-1}), $|V_p|=n-1$. We may assume that $|R \cup X_1|\geq |B \cup X_2|$.  Then $|R \cup X_1|\geq n$ because $|R \cup X_1|+ |B \cup X_2|\ge 2n$. Then for any two distinct vertices $a,  b\in R\cup X_1$, there exists a red $(a,b)$-path   on $2n-1$ vertices using edges between $R\cup X_1$ and $V_p$.  Since $(G,\tau)$ contains no red $C_{2n}$, we see that no vertex in $B $ is red-adjacent to two vertices in   $R \cup X_1$.  It follows that $|X_1|=0$ because $R$ is not blue-complete to $B$ in $(G, \tau)$. 
Suppose $ |B\cup X_2|\ge  n$. By a similar argument, $|X_2|=0$. It is easy to see that there exists a vertex in $R$ which is blue-adjacent to two vertices, say $x, y$,  in $B$. This, together with a blue $(x, y)$-path on $2n-1$ vertices using edges between $B$ and $V_p$,  yields a blue $C_{2n }$ in $(G,\tau)$.  Thus $ |B\cup X_2|\le   n-1 $ and so $|R|\ge n+1$.  Then   $(G[R],\tau) $ contains at most one red edge, else we obtain a red $C_{2n}$ using two red edges in $G[R]$ and edges between $R$ and $V_p$; and $(G[R\cup B],\tau)$ has no    red $P_3$ with both ends in $R$.   Furthermore, for any $V_\ell\in \mathcal{V}_r$ with $|V_\ell|\ge2$,  $V_\ell$ is blue-complete to $(R\less V_\ell)\cup B\cup X_2$, else $(G[R\cup B],\tau)$ contains a red $P_3$ with both ends in $R$.
Suppose $ |B\cup X_2|=  n-1\ge2$. 
 Let $x,  y\in B\cup X_2$  be  two distinct vertices. Let $P$ be a blue  $(x,y)$-path on $2n-3$ vertices using edges between $B\cup X_2$ and $V_p$.  If there exists a $V_\ell\in \mathcal{V}_r$ with $|V_\ell|\ge2$, then $V_\ell$ is blue-complete to $(R\less V_\ell)\cup B\cup X_2$. Let $u, v\in V_\ell$ and $z\in R\less V_\ell$. Then we obtain a blue $C_{2n}$ in $(G,\tau)$ with edge set  $\{zu, zv, ux, vy\}\cup E(P)$. Thus $|\mathcal{V}_r|=|R|\ge n+1\ge4$. Then   $(G[R],\tau)$  contains at most one red edge  and all other edges are colored blue. It can be easily checked that $(G[R],\tau)$ contains   a blue  $P_3$  with vertices  $u, v, w$ in order such that $ux$ and $wy$ are colored blue in $(G,\tau)$. But then $(G,\tau)$  contains a blue $C_{2n}$    with edge set  $\{uv, vw,   ux, wy\}\cup E(P)$. This proves that $ |B\cup X_2|\le  n-2 $. 
Then  $|R|\ge n+2\ge |B\cup X_2|+4$.      Let $ H $ be the   subgraph of  $G$ with $V(H)=R \cup  B \cup X_2$   and   $E(H)$  being the set of all blue edges in $(G[R \cup B \cup X_2],\tau)$. Since $G[R\cup B]$ contains no  red $P_3$ with both ends in $R$, and for any $V_\ell\in \mathcal{V}_r$ with $|V_\ell|\ge2$,  $V_\ell$ is blue-complete to $(R\less V_\ell)\cup B\cup X_2$, it follows that  for all $v\in  B\cup X_2$, $d_H(v) \ge |R|-1\ge (|H|+2)/2$;    and for all $v\in  R$, $d_H(v) \ge |H| /2$. By Theorem~\ref{Bon}, $ H $   contains a  cycle of length   $ 2n $,  which yields a blue $C_{2n}$ in $(G,\tau)$.  \qed\\

\noindent {\refstepcounter{counter}\label{n}  (\arabic{counter}) }  
$n\ge4$. 
 
 \pf Suppose $n=3$. Then $|G|\ge8$.  By  (\ref{Vpn-2}), $|V_p|=1$.  Then every edge of $(G\less X, \tau)$ is colored red or blue. Since at least one edge of $(G,\tau)$ is colored green, we may assume that  $G[X_1]$  contains a  green  edge, say $xy$.         Note that  $R$  is not blue-complete to $B$ in $(G, \tau)$ and  $|R \cup B |\geq 3$. Let $u,v,w\in R\cup B$ be all distinct such that   $u\in R$ and $v\in B$ with $uv$ colored red under $\tau$. Then we obtain  a red copy of $C_6$ in $(G, \tau)$   with vertices $x, w, y, v, u, z$ in order, where $\{z\}= V_p$, contrary to  (\ref{noC2n}).    \qed \\

\noindent {\refstepcounter{counter}\label{X1X2}  (\arabic{counter}) }  
$X_1=\emptyset $ or $X_2=\emptyset$.    

\pf Suppose   $X_1\ne\es$ and $X_2\ne\es$. Let $H:=G\less \{v_p, x_1,x_2\}$, where $v_p\in V_p$,  $x_1\in X_1$ and $x_2\in X_2$. Let $\tau'$ be obtained from $\tau$ by recoloring all the edges of  $G[V_i]$ green for all $i\in[p]$. Then $(H, \tau')$ is  Gallai $k$-colored   with no monochromatic copy of $C_{2n}$.   Since $|V_p|\le n-2$, we see that       $q_{\tau'}(H, n-1)\le q_\tau(H, n-1)=2$. Then  \[|H|=|G|-3\geq 3n-4\ge ((n-1)-1)\cdot q_{\tau'}(H, n-1)+(n-1)+1.\]   By (\ref{n}), $n-1\ge3$. By   minimality of $n$, $(H, \tau')$ contains a monochromatic, say red,  copy of $ C:=C_{2n-2}$. 
  We first claim that $C$ contains no vertex in $R\cup V_p $. Suppose   there exists a vertex $u \in R\cup V_p $ such that $u$ lies on  $C $. Let $v $ be  one neighbor  of $u$ on the cycle $C$.  Then $uv$ is colored red under $\tau'$. Since all the edges in $G[V_i]$ are colored green in $(G, \tau')$ for all $i\in[p]$, we see that either  $u \notin V_p$ or $v \notin V_p$, say the latter.  By the choice of $u$,  $v\in R\cup B\cup X_1$.  Then  we   obtain a red $C_{2n}$ in $(G,\tau)$   from the cycle $C$ by replacing $uv$ with the path having vertices $u, x_1,   v_p, v$ in order if $v\in R\cup X_1$; and with the path having vertices $u, v_p, x_1, v$ in order if $v\in B$ (and thus $u\in R$), contrary to  (\ref{noC2n}).    We next claim that  $C$ contains no vertex in $X_1$.  Suppose $V(C)\cap X_1\ne \es$.  Since $X_2$ is blue-complete to $B$ and $|X_i\less x_i|\le n-2$ for all $i\in[2]$, there must exist    $u\in   X_1\less x_1$  and  $v\in B$ such that $uv\in E(C)$.  But then we     obtain a red $C_{2n}$ in $(G, \tau)$ from the cycle $C$ by replacing $uv$ with the path having vertices $u, v_p,   x_1, v$ in order. This proves that $C$ contains no vertex in $R\cup V_p \cup X_1$.  Thus  $V(C)\subseteq B$.  By  the choice of $p$ and Lemma~\ref{k'}(ii),    $G_r^\tau\less X$  is  connected.   Thus there exist $u\in R\cup (B\less V(C))$ and $v\in V(C)$ such that $uv $ is colored red under $\tau$. Let $w$ be one neighbor of $v$ on the cycle $C$. We   obtain a red $C_{2n}$  in $(G, \tau)$ from the cycle $C$ by replacing $vw$ with the path having vertices $v, u, x_1,   w $ in order, contrary to  (\ref{noC2n}).    \qed\\

   By (\ref{X1X2}),  $X=X_1$ or $X=X_2$. 
  For the remainder of the proof,  let $\sigma$ be  obtained from $\tau$ by recoloring all the edges of $ G[X_1] $ blue if $X=X_1$,    all the edges of $ G[X_2] $ red if $X=X_2$,  and    all the edges of  $ G[V_i] $ green for all $i\in[p]$. Then $(G, \sigma)$ is   Gallai $3$-colored   with no monochromatic copy of $C_{2n}$. Note that    $q_\sigma(G, n)=q=2$ and   $\{V_1, \ldots, V_p\}$ is a Gallai partition of  $(G\less X, \sigma)$.  By (\ref{Vpn-2}),  $(G, \sigma)$ contains  no green cycle on  more than  $n-2$ vertices.  
 \bigskip
 
 \noindent {\refstepcounter{counter}\label{Vp2}  (\arabic{counter}) } 
 $  |V_p|\ge2$.
 
 \pf Suppose $  |V_p|=1$. Then $(G, \sigma)$ is Gallai $2$-colored. By Theorem~\ref{C2n},  $(G, \sigma)$  contains a red or blue cycle $C:=C_{2n}$.  By (\ref{X1X2}), we may assume that $X_2=\emptyset$. 
  By (\ref{Xi}) and the choice of $\sigma$, the cycle $C$ contains no edge  of $ G[X_1] $. It follows that $C$ is a red or blue $C_{2n}$ in  $(G, \tau)$, contrary to  (\ref{noC2n}). 
    \qed\\

 Let $x\in V_p$, $y\in R \cup X_1$ and $z\in B\cup X_2 $ be  such that $y\in X_1$ if $X_1 \ne\es$ and $z\in X_2$ if $X_2 \ne\es$. Let $H:=G\less \{x,y,z\}$.   By  (\ref{Vpn-2}) and  (\ref{X1X2}),   $q_\sigma(H, n-1)\leq q_\tau(H, n-1) =2$.    Then  \[|H|=|G|-3\geq 3n-4\ge ((n-1)-1)\cdot q_{\sigma}(H, n-1)+(n-1)+1.\]   By (\ref{n}), $n-1\ge3$. By   minimality of $n$, $(H, \sigma)$ contains a red or blue copy of $ C:=C_{2n-2}$ with vertices, say $a_1, a_2, \ldots, a_{2n-2}$ in order. We may further assume that the cycle $C$ is blue.    By (\ref{X1X2}), $C$   contains no vertex in $X_1$. 
  By the choice of $\sigma$,  the cycle  $C$   contains  no edge  in $(G[X_2], \sigma)$.    We claim that $V(C)\subseteq V_p\cup R \cup B $. Suppose  the cycle $C$  contains a vertex in $X_2$,     say  $a_1\in  X_2$.     By (\ref{X1X2}),   $X_1=\emptyset$.  By the choice of $\sigma$,  $a_2\in R \cup B \cup V_p$. By the choice of $z$, we see that $z\in X_2\less a_1$. But then we obtain a blue $C_{2n}$ in $(G, \tau)$ with vertices $a_1, y, z, a_2, \ldots, a_{2n-2}$ in order, contrary to  (\ref{noC2n}).  This proves that $V(C)\subseteq V_p\cup R \cup B $. Then $V_p\cap V(C)=\es$, else say $a_1\in V_p\less x$, then $a_2x$ is blue because $a_2a_1$ is blue, and  thus   $(G, \tau)$  has  a blue  $C_{2n}$ with vertices $a_1, z, x, a_2, \ldots, a_{2n-2}$ in order, contrary to  (\ref{noC2n}). It follows that $V(C)\subseteq    R \cup B $.   \medskip

 For the remainder of the proof, we say  a vertex $v$ on the cycle $C$ is \dfn{special} if there exist a vertex $u\in   (R\cup B) \less  V(C) $ and $V_j\in\{V_1, \ldots, V_{p-1}\}$ such that $ u,v \in V_j$, that is, $uv$ is colored green under $\sigma$.  Let $S$ be the set of all special vertices on the cycle $C$ and let $T$
 be the set of  all vertices  $u\in   (R \cup B)\less V(C)$ such that     $u$ is green-adjacent to   some $v\in S$. Finally, let $L:=V(G)\less (V(C) \cup T)$.  Then no vertex in $L$ is green-adjacent to any vertex in $V(C)\cup T$. Moreover, $V_p\subseteq L$ because $V(C)\subseteq R\cup B$.
 By the choice of $\sigma$ and the fact that $(G,\sigma)$ has neither red nor blue $C_{2n}$, we see that   no two consecutive vertices on $C$ are both special vertices.  Let $s:=|S|$. Then $s\le n-1$.  It is worth noting that for each  $V_j\in\{V_1, \ldots, V_{p-1}\}$, if $V_j\cap  S\ne \es$, then $V_j\subseteq S\cup T$ and  $|V_j|\ge2$.    Since $|V_p|\ge2$, we see that the cycle $C$ does not contain two edges  in $(G[B],\sigma)$, else we obtain a blue $C_{2n}$ in $(G, \sigma)$.  
  \\

\noindent {\refstepcounter{counter}\label{C}  (\arabic{counter}) }  
For each  $a_i\in V(C)$, if $a_i\in R$, then either $a_{i+1}\in R$ or $a_{i-1}\in R$, where all arithmetic
on indices here and henceforth is done modulo   $2n-2$.   Moreover, $|V(C)\cap R|\geq n$ and $B\less V(C)\subseteq L$.

\pf Suppose 
there exists  a vertex  $a_i \in V(C)$ such that     $a_i\in R$ but $a_{i-1},a_{i+1}\in B$.   Let $x’\in V_p\less x$. Note that $V_p\cap V(C)=\es$. But then we   obtain a  blue $C_{2n}$ in $(G, \sigma)$ with vertices $a_1,\ldots,a_{i-1},x,z,x',a_{i+1},\ldots,a_{2n-2}$ in order, contrary to  (\ref{noC2n}).  It follows that    $|V(C)\cap R|\geq n$, because   $C$ does not contain two edges  in $(G[B],\sigma)$.   Next, suppose  there exists a vertex 
   $u\in B\less V(C)$ such that $u\in T$. We may assume that   $a_1, u \in V_j$ for some $j\in[p]$. Then  we obtain  a blue $C_{2n}$ in $(G, \sigma)$ with vertices $a_1, x,u,a_2,\ldots,a_{2n-2}$ in order, contrary to  (\ref{noC2n}). Thus $B\less V(C)\subseteq L$.   \qed\\

\noindent {\refstepcounter{counter}\label{red-complete}  (\arabic{counter}) }  
Every vertex in $ L$  is red-complete to $S\cup T$ in $(G, \sigma)$,  and  for all $i,j\in[p]$ with $i\ne j$,  if $V_i  \subseteq S\cup T$  and $ V_j \subseteq S\cup T$, then    $V_i   $ is red-complete to  $ V_j  $ in $(G, \sigma)$.

\pf Suppose  there exists a vertex $v\in L$  such that $v$ is not red-complete to $S\cup T$ in $(G, \sigma)$. By the choice of $S$ and $T$, $v$ is not green-adjacent to any vertex in $S\cup T$ and thus   there must exist  some  $V_\ell\subseteq S\cup T$ such that $v$ is blue-complete to $V_\ell$. We may assume that $a_j, u \in V_\ell$ for some $j\in[2n-2]$ and    $u\in T$. But then $a_{j-1}u$ is colored blue and we obtain  a blue $C_{2n}$ in $(G,\sigma)$ with vertices $a_1,a_2,\ldots,a_{j-1},u, v,a_j,\ldots, a_{2n-2}$ in order, contrary to  (\ref{noC2n}). This proves that every vertex in $ L$  is red-complete to $S\cup T$ in $(G, \sigma)$. Next, 
suppose there exist   $V_i  \subseteq S\cup T$ and $V_j  \subseteq S\cup T$ such that $V_i$ is blue-complete to   $ V_j  $ in $(G, \sigma)$.   We may assume that  $\{a_1, u_1\}\subseteq V_i$ and $\{a_{\ell}, u_2\}\subseteq V_j$ for some $\ell\in\{2, \ldots, 2n-2\}$ and  $u_1, u_2\in T$.   But then $u_2a_{\ell+1}$ is colored blue and we obtain  a blue   $C_{2n}$ in $(G, \sigma)$ with vertices $a_1, \ldots, a_\ell, u_1, u_2,a_{\ell+1},\ldots,a_{2n-2} $ in order, 
 contrary to  (\ref{noC2n}).  \qed\\

  \noindent {\refstepcounter{counter}\label{noblue}  (\arabic{counter}) }  
If a vertex  $v\in T\cup L$ is blue-complete to an edge, say $a_1 a_{2}$, on the cycle $C$, then no vertex in $(T\cup L)\less v$ is blue-complete to any edge on the cycle with vertices $a_1,   v, a_2,\ldots, a_{2n-2}$ in order.   \bigskip

\pf It follows from the fact that    $(G, \sigma)$ has no blue $C_{2n}$.\qed\\

   For each vertex $u\in T\cup L$,  let 
\[\begin{split}
N_b(u)&:=\{v\in V(C)\mid uv \text{ is colored blue in } (G, \sigma)\} \text {  and} \\
N_r(u)&:=\{v\in V(C)\mid uv \text{ is colored red in } (G, \sigma)\}.
\end{split}
\]
 Since no vertex in $L$ is green-adjacent to any vertex on the cycle $C$, we see that   $|N_r(u)|+|N_b(u)|=2n-2$ for all $u\in L$.    By  (\ref{C}), $|N_r(u)|= |V(C)\cap R|\ge n$ for all $u\in V_p$. \\

 \noindent {\refstepcounter{counter}\label{S}  (\arabic{counter}) }  
$  |L|\le n$. Consequently,  $|T|\ge1$  and  so $s\geq 1$.

\pf Suppose for a contradiction that  $|L|\geq n+1$.  Let $L^*$ be the set of all vertices $u\in L $ with $|N_r(u)|\ge n$.    
 Let  $H^*$ be the bipartite  subgraph   of $G$  with bipartition $\{L^*, V(C) \}$ and  $E(H^*)$  consisting of all red edges between $L^*$ and $V(C)$ in $(G, \sigma)$.  By the choice of $\sigma$ and the fact that $V(C)\subseteq R\cup B$,  we see that  $H^*$ contains no cycle of length  $2n$. By Corollary \ref{C2l} applied to $H^*$ with $M=L^*$,  $N=V(C)$ and $\delta(M)\ge n>(|N|+1)/2$, we have $|L^*|\le n-1$.  It follows that $|L\less L^*|\ge 2$ and for any $v\in L\less L^*$, $|N_r(v)|\le n-1$ and so $|N_b(v)|\ge n-1$. Then $(L\less L^*)\cap V_p=\es$ by (\ref{C}).   We next claim that   there exist   two distinct vertices     $w_1, w_2\in L\less L^*$ with  $N_r(w_1)=N_r(w_2)$ such that $|N_r(w_1)|=|N_r(w_2)|=n-1$ and   every vertex in $L\less \{w_1, w_2\}$ is  red-adjacent to   at least $n-2$ vertices  in  $N_r(w_1)$. \medskip

 Suppose first that no vertex in  $  L\less L^*$ is blue-complete to any edge on   the cycle $C$. Then  for any $v\in L\less L^*$, $|N_b(v)|=n-1$ and so $|N_r(v)|=n-1$. Let   $w_1, w_2$ be two distinct  vertices in  $   L\less L^*$.  We may assume that $w_1$ is blue-complete to   $\{a_1, a_3, \ldots, a_{2n-3}\}$ in $(G, \sigma)$. Then $L\less L^*$ must be blue-complete to $\{a_1, a_3, \ldots, a_{2n-3}\}$ and   no vertex in $ L^*$ is blue-adjacent to two vertices in $\{a_2, a_4, \ldots, a_{2n-2}\}$,
 else in each case we obtain a blue $C_{2n}$ in $(G, \sigma)$.  Then $N_r(w_1)=N_r(w_2)=\{a_2, a_4, \ldots, a_{2n-2}\}$ and every vertex in $ L\less \{w_1, w_2\}$ is   red-adjacent to   at least $n-2$ vertices  in  $N_r(w_1)$, as claimed. 
 Suppose next that some vertex in      $    L\less L^*$   is   blue-complete to an edge, say  $ a_1a_2 $,  on  the cycle $C$.  Let $w_1$ be such a vertex, and  let $C^*$ be the blue cycle on $2n-1$ vertices with edge set  $\{w_1a_1, w_1a_2\}\cup  E(C \less  a_1a_2 )$. Suppose first $\{w_1\}\in \{V_1, \ldots, V_{p-1}\}$.  By  (\ref{noblue}),  every vertex in $L\less w_1$ is red-adjacent to at least $n$ vertices on $C^*$. 
If  each vertex in $V_p$ is red-adjacent to at least $n+1$ vertices on $C^*$,  let  $H^*$ be the bipartite  subgraph   of $G$  with bipartition $\{L\less w_1, V(C^*) \}$ and  $E(H^*)$  consisting of all red edges between $L\less w_1$ and $V(C^*)$ in $(G, \sigma)$. By Corollary \ref{C2l} applied to $H^*$ with $M=L\less w_1$,  $N=V(C^*)$,  $\Delta(M)\ge n+1$ and $\delta(M)\ge n =(|N|+1)/2$, $H^*$ has a cycle of length $2n$, which yields a red $C_{2n}$ in $(G, \sigma)$, contrary to  (\ref{noC2n}). Thus each vertex in $V_p$ is red-adjacent to at most   $n $ vertices on $C^*$. Since $|N_r(u)|\ge n$ for each $u\in V_p$, we see that 
$w_1$ is blue-complete to $V_p$ in $(G, \sigma)$ and $V_p$ is red-adjacent to exactly $n$ vertices on $C$.  Then  $V_p$ must be red-complete to $\{a_1, a_2\}$ in $(G, \sigma)$ and no vertex in  $ V_p$ is blue-complete to any edge on $C$, else we obtain a blue $C_{2n}$ in $(G, \sigma)$. Thus one of $a_3$ and $a_{2n-2}$, say $a_3$, must be    blue-complete to $V_p$  under  $ \sigma$. But then $(G, \sigma)$ contains a blue $C_{2n}$ with edge set $E(C\less a_2a_3)\cup \{a_2w_1, w_1x, xa_3\}$, contrary to  (\ref{noC2n}). 
This proves that $\{w_1\}\notin \{V_1, \ldots, V_{p-1}\}$. Thus there exists $V_\ell \in \{V_1, \ldots, V_{p-1}\}$ such that   $w_1\in  V_\ell$  and $|V_\ell|\ge2$. Let $w_2\in V_\ell$ with $w_2\ne w_1$. 
 By  (\ref{noblue}),   no vertex in $V_\ell$ is blue-complete to any edge on $C\less a_1a_2$.      
  This, together with fact that  $|N_b(w_i)|\ge n-1$,   implies that $|N_b(w_i)|= n-1$ for all $i\in[2]$, 
  and no vertex in $ L\less \{w_1, w_2\}$ is blue-adjacent to two vertices in $N_r(w_1)$, otherwise   we obtain  a blue   $C_{2n}$  in $(G, \sigma)$. Then    every vertex  in  $ L\less \{w_1, w_2\}$ is  red-adjacent to   at least $n-2$ vertices  in  $N_r(w_1)$,  as claimed. \medskip

  Let  $L'$   be a subset of $L$ such that $|L'|=n-1$, $ w_1, w_2\in L'$ and $(V_p\less x)\cap L'\neq \emptyset$. Then   $|L'|+|N_r(w_1)|=(n-1)+(n-1)=2n-2$.  Let $H^*$ be the bipartite subgraph of $G$ with bipartition $\{L', N_r(w_1)\}$ and $E(H^*)$ consisting of all red edges between $L'$ and $N_r(w_1)$ in $(G,\sigma)$. Then by the choice of $w_1, w_2$, we see that $\{w_1, w_2\}$ is red-complete to $N_r(w_1)$.  By Corollary~\ref{path} applied to $H^*$ with $M=L'$ and $N=N_r(w_1)$,   $H^*$ has a  cycle $C^*$ on  $2n-2$ vertices. By the choice of $L'$, let $x_1x_2$ be an edge of $C^*$ with $x_1 \in  V_p\less x $. Let   $u\in  V(C)\cap R$ with $u\notin V(C^*)$. But then   $(G, \sigma)$  has   a red $C_{2n}$  with edge set $E(C^*\less x_1x_2 )\cup \{x_1u,ux,xx_2\}$, contrary to  (\ref{noC2n}). 
This proves that $|L|\le n$. Consequently, $|T|=|G|-|V(C)\cup L|\ge 3n-1-(2n-2+n)\ge1$ and so $s=|S|\ge1$.\qed \\

By (\ref{red-complete}) and  (\ref{S}), $X_2=\es$.   By (\ref{C}),    $z\in L\less V_p$. We next claim that \\

\noindent {\refstepcounter{counter}\label{s}  (\arabic{counter}) }  
$s\le n-4$.

 \pf Suppose  that $s\ge n-3$.    
  By  (\ref{red-complete}) and the fact that $V_p\subseteq L$ and $|V_p|\le n-2$, we may assume that $S\cup T=V_1\cup \cdots\cup V_\ell$ for some $ \ell\in[p-1]$.   By  (\ref{red-complete}) again, $L$ is red-complete to $V_1\cup \cdots\cup V_\ell$ and $V_i$ is red-complete to $V_j$ for all $i, j\in[\ell]$ with $i\ne j$ when $\ell\ge2$.   By  (\ref{C}), $|V(C)\cap R|\ge n$ and so $(V(C)\less S)\cap R\ne \es$.    Let $X^*:=\{x\}$ when $s\ge n-2$ and $X^*:=\{x, x'\}$ when $s=n-3$, where      $x'\in V_p\less x$.   Then   
 \[
 \begin{split}
 |S\cup T\cup (L\less X^*)|&=|G|-|X^*\cup (V(C)\less S)|\\
 &\ge 
 \begin{cases}  
(3n-1)-(n+1)=2n-2   &  \text{if } |X^*|=1 \\
(3n-1)-(n+3) =2n-4   &  \text{if } |X^*|=2\\
\end{cases}\\	
 &=2n-2|X^*|.
 \end{split}\]  
 We define  $L^*$  to be a subset of $S\cup T\cup (L\less X^*)$ with  $|L^*|=2n-2|X^*|$ such that when $|X^*|=1$,  $L^*\cap V_p\ne \es$; and  when   $|X^*|=2$,    $z\in L^* $ and $L^*\cap T\ne\es$.    Let $H^* $ be the subgraph of $G$ with   $V(H^*)=L^*$ and $E(H^*)$ consisting of   all red edges in $(G, \sigma)$ between $L^*\cap (L\less X^*)$   and $L^*\cap(V_1\cup \cdots\cup V_\ell)$,   and all  edges between each pair $L^*\cap V_i$ and $L^*\cap V_j$ for all $i, j\in[\ell]$ with $i\ne j$.  Since $|V_i|\le n-2$ for all $i\in[\ell]$,  $|L\less X^*|\le n-|X^*|$, we see that       $\delta(H^*)\ge |H^*|/ 2 $. When   $X^*=\{x\}$, let $u\in  (V(C)\cap R)\less S $.   By Theorem \ref{Bon}, $H^*$ contains a   cycle   $C^*:=C_{2n-2}$, say  with vertices $v_1,v_2, \ldots,v_{2n-2}$ in order. We may assume that $v_1\in V_p$. Then $\{u, v_2\}$ is red-complete to $\{x,v_1\}\subseteq V_p$ and $(G, \sigma)$ contains  a red $C_{2n}$   with vertices $v_1, u, x, v_2, \ldots, v_{2n-2}$ in order, contrary to  (\ref{noC2n}). Thus  $X^*=\{x, x'\}$ and $s=n-3$.  By Theorem \ref{Bon} again, $H^*$ contains a   cycle   $C^*:=C_{2n-4}$, say  with vertices $v_1,v_2, \ldots,v_{2n-4}$ in order. We may further assume that $v_1=z$. By the choice of $E(H^*)$,   $v_2\in S\cup T$. Note that $v_1\in B$ by the choice of $z$.  
\medskip

Since  $C\less S$ contains  at most $s-1$ isolated vertices and $2n-2-(s+s-1)= 2n-1-2(n-3)=5$, it follows that  $C\less S$ contains at least three edges.   By (\ref{noblue}), $V_p$ is not blue-complete to any two edges on the cycle $C$. Let  $ab, cd\in E(C\less S) $ such that $a, c\in R$. 
 Then $b, d\in R$, else,  say $b\in B$, then $(G, \sigma)$ contains a blue $C_{2n}$ with edge set $E(C\less ab)\cup\{av_1, v_1x, xb\}$ if $av_1$ is blue; and a red $C_{2n}$  with vertices $v_1, a, x, c, x', v_2, \ldots, v_{2n-4}$ in order   if $av_1$ is red, contrary to  (\ref{noC2n}).  Thus   $v_1$ is blue-complete to $\{a,b,c,d\}$, else,  say $v_1a$ is red,     then $(G, \sigma)$ contains  a red    $C_{2n}$ with vertices $v_1,a,x,c,x',v_2,\ldots,v_{2n-4}$ in order, contrary to  (\ref{noC2n}).  
  By (\ref{noblue}), no vertex in $T\cup (L\less z)$  is blue-complete to $\{a,b\}$ or  $\{c,d\}$ in $(G, \sigma)$.  
 By the choice of $S\cup T$, we see that no vertex in $S$ is blue-complete to $\{a,b\}$ or  $\{c,d\}$ in $(G, \sigma)$. We may assume that $v_3a$ is colored red because   $v_3\in S\cup T\cup (L\less z)$. 
 But then   $(G, \sigma)$ contains  a red    $C_{2n}$ with vertices $v_1,  v_2,  x, b, x', a, v_3,    \ldots,v_{2n-4}$ in order, contrary to  (\ref{noC2n}). \qed\\

 \noindent {\refstepcounter{counter}\label{S*}  (\arabic{counter}) }  
There exists  a subgraph $J$ of $C\less S$    with  $| J | \ge 2n-2-2s\ge6$ such that each component of $J$ is a path of even order, and   for all $u\in   T\cup L$ except possiblely one vertex, $u$ is neither  blue-complete to an  edge in  $J$ nor green-adjacent to a  vertex in $V(J)$, and so    
    $|N_r(u)\cap  V(J)|\ge |J|/2$.     \bigskip
 
\pf  By  (\ref{S}), $s\ge1$ and $|T|\ge1$. 
 By  (\ref{s}), $|V(C)\less S|\ge n+2$.   Since    $S$     is an independent set of the cycle $C$,       
we see that  $C\less S$ has exactly $s$ components such that    each component is a path. Assume first that  no vertex in     $T\cup L $ is blue-complete  to an edge on the cycle  $C $. Let  $J$ be  obtained from  $C\less S$ by  deleting one end of each odd component of  $C\less S$. Then $| J|\ge 2n-2-2s $ and for all $u\in T\cup L$, $u$ is not green-adjacent  to  any vertex in  $ V(J)$,  and so $|N_r(u)\cap  V(J)|\ge | J |/2$ because $u$ is not blue-complete  to any edge on the cycle  $C $.  By   (\ref{s}),   $| J|\ge 2n-2-2s\ge 6$.
Assume next that there exists a vertex $w\in T\cup L$ such that $w$ is blue-complete to an edge, say $a_1a_2$, on the cycle $C$.  By (\ref{noblue}),   no vertex in $(T\cup L)\less w$ is blue-complete to any edge  on the path $P:=C\less a_1a_2$.  Let   $J$ be   obtained from  $P\less S$ by  deleting one end of each odd component of  $P\less S$.   We claim that $P\less S$ has at most  $s$ odd components.  Suppose $P\less S$ has at least  $s+1$ odd components. Then $P\less S$ has exactly   $s+1$ components such that each component is odd. But then $P  $ must have an odd number of vertices, contrary to $|P |=2n-2$.   Thus $P\less S$ has at most  $s$ odd components, as claimed.  Then   $| J|\ge 2n-2-2s\ge6$ and  for all $u\in  T\cup L$ with $u\ne w$, $u$ is not green-adjacent  to  any vertex in  $ V(J)$, and so   $|N_r(u)\cap  V(J)|\ge |J|/2$  because $u$ is not blue-complete to any edge in $J$. \qed\\
 
  By  (\ref{S*}), let $w$ be the possible vertex in $T\cup L$ such that $w$ is blue-complete to some edge in $J$.     Let $T':=T\less w $ and $L':=L\less w $. Then $|T'\cup L'|=|(T\cup L)\less w|\ge n$.  We next show that $T'\ne \es$. Suppose $T'= \es$. Then $T=\{w\}$ and $L=L'$ with $|L|\ge n$.      By  (\ref{S*}),  for any $u\in L$,  $|N_r(u)\cap  V(J)|\ge |J|/2 \ge n-1-s$ and $u$ is not green-adjacent to any vertex in $V(J)$.  
  By  (\ref{red-complete}),    $L$ is red-complete to $S\cup T$.  Let $H^*$ be the bipartite graph of $G$ with bipartition $\{L, V(C)\cup \{w\}\}$ and $E(H^*)$ consisting of all red edges between $L$ and $V(C)\cup \{w\}$ in $(G, \sigma)$.  Then each vertex in $L\less V_p$ is red-adjacent to at least $|N_r(u)\cap  V(J)|+|S\cup T|\ge (n-1-s)+s+1=n$ vertices in $V(C)\cup \{w\}$, and each vertex in $V_p$ is red-adjacent to at least $|T|+|R\cap V(C)|\ge n+1$ vertices in $V(C)\cup \{w\}$. By Corollary \ref{C2l} applied to $H^*$ with $M=L$, $N=V(C)\cup\{w\}$,  $\delta(M)\ge n$ and $\Delta(M)\ge n+1$, $H^*$ has a cycle of length $2n$, which yields a red $C_{2n}$ in $(G, \sigma)$, contrary to  (\ref{noC2n}).  Thus  $T'\ne \es$.  \medskip

  Let $x_1\in V_p $ with $x_1\ne x$. Note that   $|T'\cup L' |\ge n= (n-s-1)+ (s+1)\ge (n-s-1)+2$.   Let $W$ be a subset of  $T'\cup L'$ with  $|W|=n-s-1$ such that $ x, x_1\notin W$,  $|W\cap T'|$ is  as large as possible and $|W\cap V_p|$ is as small as possible.   Let $W^*$ be a subset of $(T'\cup L')\less W$ with $|W^*|=s$ such that $x_1\in W^*$, $x\notin W^*$  and $|W^*\cap V_p|$ is as large as possible. 
 Let $H^*$ be the subgraph of $G$ with  $V(H^*)= W^*\cup S $ and $E(H^*)$ consisting of all red edges   in $(G[W^*\cup S], \sigma)$.   We claim that $\delta(H^*)\ge|H^*|/2$. Assume first $|V_p |\ge s+1 $. By the choice of $W$ and $W^*$, $W^*\subseteq V_p$.  Then  $W^*$ is red-complete to $S$ and so $\delta(H^*)\ge|H^*|/2$. Assume next $|V_p |\le s $.  By  (\ref{red-complete}), every vertex in $W^*\cap L'$ is red-complete to $S$;   for  every vertex   $u\in W^*\cap T'$, we may assume that  $u\in V_\ell$ for some $\ell\in[p-1]$. Then $u$ is red-complete to $V(H^*)\less V_\ell$, and so   $\delta(H^*)\ge|H^*|/2$ because   $|V_\ell|\le |V_p    |\le s  $. Hence in both cases,    $\delta(H^*)\ge|H^*|/2$.  By Theorem \ref{Bon}, $H^*$  has  a cycle $C^*$ on $2s$ vertices (here by abusing the notation, $C^*$ denotes an edge on two vertices when $s=1$). Next, 
 let $H'$ be the bipartite subgraph of $G$ with bipartition $\{W, V(J)\}$ and $E(H')$ consisting of all red edges between $W$ and $V(J)$ in $(G, \sigma)$.   
  By  (\ref{S*}),   every vertex in $W$ is red-adjacent to at least $ |V(J)| /2$ vertices in $V(J)$. Suppose first that $W$ is a disjoint union of $M_1\ne\es$ and $M_2\ne\es $ and $V(J)$ is a disjoint union of $N_1$ and $N_2$ such that 
  $|N_1|=|N_2|=|V(J)|/2 \ge n-1-s\ge3$ and  $M_i  $ is red-complete to $N_i$ but blue-complete to $N_{3-i}$ for all $i\in[2]$.   
   By  (\ref{S*}) and the choice of $W$,   neither $J[N_1]$ nor $J[N_2]$ contains an edge  and  $ J$ contains at least three  independent edges, say $a_ia_{i+1}, a_ja_{j+1},a_\ell a_{\ell+1}$ with $i<j<\ell$.  We may further assume $a_i,a_j\in N_1$. Then $a_{i+1},a_{j+1}\in N_2$    and $(G, \sigma)$ has a blue $C_{2n}$ with edge set $E(C \less\{a_ia_{i+1}, a_ja_{j+1}\})\cup\{a_ib_2,b_2a_j,a_{i+1}b_1,b_1a_{j+1}\}$,  where $b_1\in M_1$ and $b_2\in M_2$, contrary to  (\ref{noC2n}).  This proves that $W$ and $V(J)$ have no such partition.  
By Corollary \ref{P2l-13}  applied to $H'$ with $M=W$, $N=V(J)$,  $\delta(M)\ge |V(J)|/2$, $H'$ has an $(a, b)$-path $P$ on $2(n-1-s)-1=2n-(2s+3)$ vertices  with $a, b\in W$.  Note that $|(V(P)\cup V(C^*))\cap V(C)|= (n-2-s)+s=n-2$.  \bigskip

 \noindent {\refstepcounter{counter}\label{noC'}  (\arabic{counter}) }  If  $(G,\sigma)$ contains   a red cycle $C_*$  on  $2n-2$ vertices  such that $x_1\in V(C_*)$ and $x\notin V(C_*)$, then $|V(C_*)\cap V(C)|\ge n$. 
 
 \pf  Suppose $|V(C_*) \cap V(C)|\le n-1$.   By  (\ref{C}),  $|R\less V(C_*)|\ge |  (R\cap V(C))\less V(C_*)|\ge 1$. Let   $u\in   R \less V(C_*) $ and $x_2$ be a neighbor of $x_1$  on $C_*$. Then $xx_2$ is red because $x,x_1\in V_p$. But then $(G, \sigma)$ has a red $C_{2n}$ with edge set $ E(C_*\less x_1x_2)\cup\{ x_1u, ux, xx_2\}$, contrary to  (\ref{noC2n}).  \qed\\

Let $\ell\in[p-1]$ such that  $V_\ell \cap T'\ne \emptyset$. This is possible  because $T'\ne \emptyset$. By the choice of $T$, $S \cap V_\ell\ne \emptyset$. For the remainder of the proof, let $b_1\in S \cap V_\ell$.   By the choice of $H^*$, $b_1, x_1\in V(C^*)$.     Let $b_2$ be a neighbor of $b_1$ on the cycle   $C^*$;   let $E_0=\{b_1b_2\}$ when $s\ge2$ and $E_0=\emptyset$ when $s=1$. 
    Since $|V(P)\cap V(J)|  =n-2-s$ and every  minimum vertex cover of $J$ has at least $n-1-s$ vertices, we see that $J\less V(P) $ has at least one edge,     say $d_1d_2$.  We next claim that  we can choose $P$ so that either $a\in T'$ or $b\in T'$. Suppose not. Then $H'$ has no cycle on $2(n-1-s)$ vertices because $T'\cap W\ne \es$. By the choice of $W$, $a, b\in L'$.   By (\ref{red-complete}), $ab_1$ and $bb_1$ are red because $b_1\in S$. Then neither $d_1$ nor $d_2$ is red-complete to $\{a,b\}$ because $H'$ has no cycle on $2(n-1-s)$ vertices.  By  (\ref{S*}) and the choice of $b_1$ and $J$,  no vertex in $ T'\cup L'$ is blue-complete to $\{d_1, d_2\}$. Thus we may assume that $ad_1, bd_2, b_1d_1$ are red.  If $b_2\in S $, then  $bb_2$ is red. But then  $(G, \sigma)$ has a red $C_*$ on $2n-2$ vertices with edge set $E(P)\cup (E(C^*)\less E_0)\cup\{b_1d_1,d_1a, bb_2\}$ such that $|V(C_*)\cap V(C) |\le |(V(P)\cup V(C^*))\cap V(C)|+1= n-1$,  contrary to (\ref{noC'}).  Thus $b_2\notin S$ and so $b_2$ is not blue-complete to $\{d_1,d_2\}$. We may  assume that $b_2d_2$ is red. But then $(G, \sigma)$ has a red $C_*$ on $2n-2$ vertices with edge set $E(P)\cup (E(C^*)\less E_0)\cup\{b_1a, bd_2, d_2b_2\}$ such that $|V(C_*)\cap V(C)|\le |(V(P)\cup V(C^*))\cap V(C)|+1= n-1$,  contrary to (\ref{noC'}).  This proves that we can choose $P$ so that either $a\in T'$ or $b\in T'$, say the former.   \medskip

By the choice of $T$ and $b_1$, we may further assume that $a\in V_\ell$ and so $ab_1$ is colored green.  Then  $ab_2$ is red.   Recall that  $a\in T'$. By  (\ref{S*}), $a$ is not blue-complete to $\{d_1, d_2\}$ and neither $ad_1$ nor $ad_2$  is green. We may assume that $ad_1$ is red. Then $b_1d_1$ is red.   
It follows that   $bb_2$ is not colored red, otherwise $(G, \sigma)$ has a red $C_*$ on $2n-2$ vertices with edge set $E(P)\cup (E(C^*)\less E_0)\cup\{b_1d_1, d_1a,  bb_2\}$ such that $|V(C_*)\cap V(C)|\le |(V(P)\cup V(C^*))\cap V(C)|+1= n-1$, contrary to (\ref{noC'}).     Similarly,  $bd_1$ is blue because $ab_2$ is red. Note that  $ad_1$ is red and $bd_1$ is blue,  we see that  $b\notin V_\ell$. Then  $b$ is red-complete to $\{a, b_1\}$  by  (\ref{red-complete}) and $bd_2$ is red  by  (\ref{S*}).  If   $b_2d_2$ is blue, then $bb_2$ is not colored green because $bd_2$ is red. Thus   $bb_2$ is   blue because $bb_2$ is not red. But then  $(G, \sigma)$ has a blue $C_{2n}$ with edge set $E(C\less d_1d_2)\cup \{d_1b, bb_2, b_2d_2\}$, contrary to  (\ref{noC2n}).  Thus $b_2d_2$ is red. Since $a\in T'$ and $b$ is red-complete to $\{a, b_1\}$, by the choice of $H'$, there must exist three consecutive vertices, say $ v_1, v_2, v_3$, on $P$  
with  $v_1, v_3\in W$ and $v_2\in V(J)$ such that  either $\{v_1, v_3\}\cap T\ne \es$ and $\{v_1, v_3\}\cap L\ne \es$,  or $v_1\in V_\ell'\cap T$ and  $v_3\in V_{\ell''}\cap T$  for some $\ell', \ell''\in [p-1]$ with $\ell'\ne \ell''$.   
By (\ref{red-complete}), $v_1v_3$ is colored red. 
But then  $(G, \sigma)$ has a red $C_*$ on $2n-2$ vertices with edge set $E(P\less v_2)\cup (E(C^*)\less E_0)\cup\{v_1v_3,  b_1d_1, d_1a, bd_2, d_2b_2\}$ such that $|V(C_*)\cap V(C)|\le |(V(P\less v_2)\cup V(C^*))\cap V(C)|+2= n-1$,  contrary to (\ref{noC'}).      \medskip

This completes   the proof of Theorem \ref{even}.
\qed

\section*{Acknowledgements}

The authors would like to thank  Christian Bosse and Jingmei Zhang for their helpful discussion.

\end{document}